\documentclass[a4paper,10pt, dvipsnames]{amsart}
\usepackage{amssymb}
\usepackage{amsthm}
\usepackage[english]{babel}
\usepackage[shortlabels]{enumitem}

\theoremstyle{plain}
\newtheorem{thm}{Theorem}[section]
\newtheorem{prop}[thm]{Proposition}
\newtheorem{lemma}[thm]{Lemma}
\newtheorem{cor}[thm]{Corollary}
\newtheorem{mainthm}{Main Theorem}
\theoremstyle{definition}
\newtheorem{alg}[thm]{Algorithm}

\newtheorem{remark}[thm]{Remark}
\newtheorem{example}[thm]{Example}

\DeclareMathOperator{\End}{End}
\DeclareMathOperator{\Hom}{Hom}

\DeclareMathOperator{\Pic}{Pic}
\DeclareMathOperator{\Mat}{Mat}
\DeclareMathOperator{\GL}{GL}

\def\Q{\mathbb{Q}}
\def\Z{\mathbb{Z}}
\def\F{\mathbb{F}}

\def\p{\mathfrak{p}}
\def\frf{\mathfrak{f}}

\newcommand{\cB}{\mathcal{B}}
\newcommand{\cI}{\mathcal{I}}
\newcommand{\cL}{\mathcal{L}}
\newcommand{\cM}{\mathcal{M}}
\newcommand{\cO}{\mathcal{O}}
\newcommand{\cQ}{\mathcal{Q}}

\newcommand{\vphi}{\varphi}
\newcommand{\set}[1]{\left\lbrace#1\right\rbrace }

\usepackage[hyphens]{url}
\usepackage{hyperref}
\usepackage[hyphenbreaks]{breakurl}

\usepackage{color}

\title[Modules over orders]{Modules over orders,\\ conjugacy classes of integral matrices, and abelian varieties over finite fields}
\author{Stefano Marseglia}
\address{{Mathematical Institute, Utrecht University, P.O. Box 80010, 3508 TA, Utrecht, The Netherlands}}
\address[current]{{Laboratoire de math\'ematiques GAATI, Université de la Polynésie Française, Faaa, French Polynesia}}
\email{{stefano.marseglia@upf.pf}}
\date{\today}


\begin{document}
\begin{abstract}
    We give an algorithm to compute representatives of the conjugacy classes of semisimple square integral matrices with given minimal and characteristic polynomials.
    We also give an algorithm to compute the~$\F_q$-isomorphism  classes of abelian varieties over a finite field~$\F_q$ which belong to an isogeny class determined by a characteristic polynomial~$h$ of Frobenius when~$h$ is ordinary, or~$q$ is prime and~$h$ has no real roots.
\end{abstract}

\maketitle

\section{Introduction}
The problem of classifying finitely generated modules over a commutative unitary ring~$R$ up to isomorphism is in general very hard.
It has been studied for several classes of rings~$R$, producing a vast literature which goes back to Eichler and Jacobinski, and spans several decades.
Excellent accounts can be found, for example, in the books~\cite{CurtisReiner81} and~\cite{CurtisReiner87} by Curtis and Reiner, in the books~\cite{RoggenkampHuber-Dyson70} by Roggenkamp and Huber-Dyson, and in~\cite{Roggenkamp70} by Roggenkamp.
Further references are provided in this introduction.

In some very special cases, we do have theoretical classifications of the isomorphism classes of modules: this is the case, for example, when~$R$ is a field, a principal ideal domain, or a Dedekind domain.
The latter case, which is usually referred to as Steinitz Theory~\cite{Steinitz12}, will be the starting point of the investigation contained in this paper.

In other cases, as explained in Section~\ref{subsec:comparison}, there are algorithms to produce representatives of the isomorphism classes.
In this paper, we produce a new algorithm that vastly supersedes the previous ones in terms of generality.

\subsection{First main result}\
In this introduction, in order to simplify the presentation, we will work under more restrictive hypotheses than in most of the paper.
We will give only the essential definitions and refer the reader to Section~\ref{sec:notation} for the missing ones.
The following is the first main result of this paper, which can be found as Theorem~\ref{thm:algcorrect} in Section~\ref{sec:modules}.
\begin{mainthm}\label{mainthm:1}
    Let~$R$ be a~$\Z$-order in a direct product of number fields $K=K_1\times\cdots\times K_n$.
    For positive integers~$s_1,\ldots,s_n$, consider the $K$-module $V=K_1^{s_1}\oplus\cdots\oplus K_n^{s_n}$.
    We provide Algorithm~\ref{alg:isom_classes} that computes representatives of the~$R$-linear isomorphism classes of~$\Z$-lattices of full rank in~$V$ that are closed under the induced action by~$R$.
\end{mainthm}
Very briefly, Algorithm~\ref{alg:isom_classes} uses the fact that all representatives of the sought isomorphism classes can be found among the~$R$-modules between two specific ones, that depend on the maximal order~$\cO$ of~$K$ and the conductor~$\frf=(R:\cO)$ of~$R$ in~$\cO$; see Theorem~\ref{thm:bwcondandmaxord}.
Theorem~\ref{thm:bwcondandmaxord} can be interpreted as an effective version of the Jordan-Zassenhaus Theorem, specialized to commutative $\Z$-orders, which states that the number of isomorphism classes is finite.
Once we have a list of these intermediate modules, one can sieve out a minimal set of representatives of the isomorphism classes using the work by Bley-Hofmann-Johnston~\cite{2022BleyHofmannJohnston}, whose relevance for our work is discussed in Remark~\ref{rmk:is_isom}.
In Remark~\ref{running_time}, we discuss the running time of Algorithm~\ref{alg:isom_classes}.
As a by-product of Algorithm~\ref{alg:isom_classes}, we obtain a scheme to sort, and hence label, these modules; see Remark~\ref{rmk:labelling}.
From our methodology it is also easy to see that the number of the isomorphism classes computed in Main Theorem~\ref{mainthm:1} is divisible by the size of~$\Pic(\cO)$. 
See~Corollary~\ref{cor:divbyclassnum}.

We pause to stress that the two facts just stated above, namely, that the classes have representatives in the range mentioned, and that the size of~$\Pic(\cO)$ divides the number of the classes, are certainly known to the experts, and can be deduced from known results about genera of modules.
Recall that the genus of a module~$M$ is the set of modules which are locally isomorphic to~$M$ at every rational prime~$p$; see for example~\cite[Sec.~31]{CurtisReiner81}.
In particular, the first statement can be deduced from~\cite[Satz~2]{Jacobinski68}; see also~\cite[Ch.~VII, Thm.~1.11]{Roggenkamp70}.
The second statement follows from the fact that the stable isomorphism classes of modules in a given genus form an abelian group, which admits a surjective map onto~$\Pic(\cO)$; see~\cite[Lemma~2.9]{Wiegand84}.
By the Jacobinski Cancellation Theorem, stable isomorphism is the same as isomorphism as long as the order satisfies the Eichler condition; see~\cite[Sec.~51]{CurtisReiner87}.
We also mention that these results have been generalized for other classes of rings and modules in a series of papers by Guralnick; see~\cite{Guralnick84},~\cite{Guralnick86},~\cite{Guralnick87}. 
Our proofs of Main Theorem~\ref{mainthm:1} and Corollary~\ref{cor:divbyclassnum} are shorter and simpler. In particular, we don't need the concepts of genera or stable isomorphism.

\subsection{Applications}\

In this paper we treat two applications of Algorithm~\ref{alg:isom_classes}.
The first one, which is treated in Section~\ref{sec:matrices}, is a solution to the conjugacy classes problem for semisimple~$r\times r$ integral matrices.
Recall that a matrix is semisimple if its minimal polynomial is squarefree, that is, with no repeated irreducible factors.
Two such matrices~$A$ and~$B$ are conjugate if there exists~$P\in \GL_r(\Z)$ such that~$PA=BP$. 
The problem of understanding the conjugacy classes with given minimal and characteristic polynomials has a long history. 
The first results are due to Latimer-MacDuffee~\cite{LaClMD33} for matrices with squarefree characteristic polynomial.
See also~\cite{taussky49}.
These two papers describe bijections between certain modules over orders and conjugacy classes of integral square matrices with prescribed characteristic polynomial.
Generalization of the Latimer-MacDuffee correspondence are provided by Bender in~\cite{Bender}, Buccino~in \cite{Buccino69}, Estes and Guralnick in~\cite{EstesGuralnick84} and Brzezinski in~\cite{BRZEZINSKI199021}.
Further results on conjugacy of matrices can be found in \cite{Rehm77}, \cite{GURALNICK1980241}, \cite{Wallace84}, \cite{Solomon99} and~\cite{Knight_Stasinski_2024}.

In~\cite[Thm~8.1]{MarICM18} we give a direct proof of a version of the Latimer-MacDuffee Theorem that applies to all square semisimple integral matrices. 
This generalization, which is recorded in the text as Theorem~\ref{thm:LM}, when combined with Algorithm~\ref{alg:isom_classes}, leads to the second main theorem of the paper, which can be found below as Theorem~\ref{thm:conjclasses}.

\begin{mainthm}\label{mainthm:2}
    Let~$m_1,\ldots,m_n$ be irreducible monic polynomials in~$\Z[x]$ generating pairwise coprime ideals in~$\Q[x]$. 
    Put~$m=m_1\cdot\ldots\cdot m_n$ and~$h=m_1^{s_1}\cdots m_n^{s_n}$, where~$s_1,\ldots,s_n$ are positive integers.  
    Algorithm~\ref{alg:isom_classes} allows us to compute the conjugacy classes of integral matrices with minimal polynomial~$m$ and characteristic polynomial~$h$.
\end{mainthm}

The second application of Algorithm~\ref{alg:isom_classes} regards abelian varieties over a finite field~$\F_q$.
It constitutes the content of Section~\ref{sec:AV}, and we refer the reader to that for missing definitions.
Consider an~$\F_q$-isogeny class~$\cI_h$ of abelian varieties over~$\F_q$ determined by a characteristic polynomial~$h\in\Z[x]$ of Frobenius.
We assume that~$\cI_h$ is ordinary, or that~$q$ is prime and~$h$ has no real roots.
Observe that we make no assumptions on the factorization of~$h$.
Using results by Deligne~\cite{Del69} and Centeleghe-Stix~\cite{CentelegheStix15}, we can describe these abelian varieties in terms of modules over a certain order determined by~$h$; see Theorem~\ref{thm:eqcatAV}.
This result, combined with Algorithm~\ref{alg:isom_classes}, leads to the third main theorem of the paper, which is later reported as Theorem~\ref{thm:isomclassesAV}.
\begin{mainthm}\label{mainthm:3}
    Let~$\cI_h$ be the~$\F_q$-isogeny class of abelian varieties over the finite field~$\F_q$ determined by the characteristic polynomial~$h$ of Frobenius.
    Assume that~$h$ is ordinary, or that~$q$ is prime and~$h$ has no real roots.   
    Algorithm~\ref{alg:isom_classes} allows us to compute the~$\F_q$-isomorphism classes of abelian varieties in~$\cI_h$.
\end{mainthm}
The above-mentioned labelling scheme for modules induces a way to sort and label the abelian varieties we are computing in Main Theorem~\ref{mainthm:3}.
This is interesting if one wants to incorporate data about the isomorphism classes of such abelian varieties in a database, like the LMFDB~\cite{lmfdb}.

Finally, in Section~\ref{sec:example} we include two examples.
In Example~\ref{ex:AV} we compute isomorphism classes of abelian varieties over~$\F_3$ in a given isogeny class.
In particular, we exhibit~$3$ pairwise non-isomorphic abelian surfaces which become all isomorphic after taking the product with one elliptic curve.
In Example~\ref{ex:matrices} we compute conjugacy classes of matrices.

The implementation in \texttt{MAGMA}~\cite{Magma} of Algorithm~\ref{alg:isom_classes} is available at 
\[ \text{\url{https://github.com/stmar89/AlgEt}}\footnote{at the moment of submission, the most recent commit is 4c22349},\]
together with the code to reproduce the examples (see the webpage of the author for a precise link).
Part of this implementation uses functionalities from the \texttt{julia} package \texttt{Nemo/Hecke} \cite{nemo}.

\subsection{Comparison with previous and related results}\
\label{subsec:comparison}
Main Theorem~\ref{mainthm:1} is a generalization of the results contained in~\cite{MarICM18}, where we consider the case~$V=K$, that is,~$s_i=1$ for all~$i$. In this case, the~$R$-modules we are considering are fractional~$R$-ideals.
In~\cite{MarICM18}, we first study the problem locally, introducing the notion of weak equivalence between fractional ideals. 
Two fractional~$R$-ideals~$I$ and~$J$ are weakly equivalent if~$I_\p\simeq J_\p$ for every prime ideal~$\p$ of~$R$.
As observed in~\cite[Sec.~5]{LevyWiegand85}, two fractional ideals are weakly equivalent if and only if they belong to the same genus.
Once we have computed the weak equivalence classes, for every overorder~$S$ of~$R$, we let the Picard group~$\Pic(S)$ act on the weak equivalence classes with multiplicator ring~$S$ to compute all the isomorphism classes of fractional~$R$-ideals.
Computing weak equivalence classes in general requires an expensive enumeration step, analogous to the one needed in Algorithm~\ref{alg:isom_classes}.
On the other hand, since weak equivalence is coarser than isomorphism, the quotient in which enumeration takes places is typically smaller.
Furthermore, checking whether two ideals are weakly equivalent is faster than checking whether they are isomorphic.
Moreover, in certain cases, we can skip the enumeration part entirely.
Indeed, if~$S$ is an overorder of~$R$ which is Gorenstein, that is, every fractional~$R$-ideal with multiplicator ring~$S$ is invertible in~$S$, then there is only one weak equivalence class with multiplicator ring~$S$.
We extend this statement in~\cite[Sec.~6]{2022MarCMtype} where we classify all weak equivalence classes with multiplicator rings which are close-to-being Gorenstein.

Another case where we have a method to efficiently compute the~$R$-modules is when~$V=K^s$ for some positive integer~$s$, that is, when~$s_i=s$ for all~$i$, and the order~$R$ is assumed to be Bass, which means that all overorders are Gorenstein.
In this case there is a classification of the~$R$-modules in~$V$ due to Bass~\cite{ba62} and Borevi\v{c}-Faddeev~\cite{BF65rus}.
Such a classification makes them easy to compute.
If~$R$ is Bass but~$V$ is not a pure power of~$K$ then the classification of the classes becomes immediately much more complicated; see~\cite{BaethSmerting22}.

In the cases~$V=K$ and~$V=K^s$ with~$R$ Bass, we used the above-mentioned results to produce algorithms to compute conjugacy classes of~$\Z$-matrices, see~\cite[Cor.~8.2]{MarICM18}, and isomorphism classes of abelian varieties over finite fields, see~\cite{MarAbVar18} and~\cite{MarBassPow}.
For the reasons we just explained, these specialized algorithms will typically perform better than the methods described in Main Theorems~\ref{mainthm:2} and~\ref{mainthm:3}. 

We conclude the introduction with two asides.
Contrary to the problem of computing representatives of the conjugacy classes of square integral matrices which is discussed in this paper, the problem of determining whether two such matrices are conjugate has recently received a lot of attention.
Algorithms can be found in the thesis of Husert~\cite{HusertPhDThesis17}, 
in previous work of the author~\cite{MarICM18}, 
in Eick-Hofmann-O'Brien~\cite{EickHofmannOBrien19},
and in the above-mentioned paper by Bley-Hofmann-Johnston~\cite{2022BleyHofmannJohnston}.
The first two deal with special cases, while third and the fourth work with all square integral matrices, including the non-semisimple ones.

The second comment regards almost-ordinary abelian varietietes over finite fields of odd characteristic. For simple isogeny classes, the work of Oswal-Shankar~\cite{OswalShankar20} gives a description of such abelian varieties in terms of modules analogous to the ones used to obtain Main Theorem~\ref{mainthm:3}. 
These results were generalized by Bergstr\"om-Karemaker and the author in~\cite[Thm.~2.12]{BergKarMarIMRNOnline} to isogeny classes of almost-ordinary abelian varieties in odd characteristic, with commutative endomorphism algebras, or, equivalently, with squarefree characteristic polynomial. 
For this reason, one should use the specialized algorithms contained in~\cite{MarICM18} and~\cite{MarAbVar18} to compute them.

\subsection*{Acknowledgements}\
The idea of this paper took form during a collaboration with Jonas Bergstr\"om and Valentijn Karemaker.
The author is thankful for their encouragement to write it down.
The author expresses his gratitude to them and to Tommy Hofmann for comments on a preliminary version.
Special thanks go to Robert Guralnick for suggesting improvements and several references. 
The author is grateful to the anonymous referees for carefully reading the paper and suggesting several improvements.
The author was supported by Nederlandse Organisatie voor Wetenschappelijk Onderzoek, grant number VI.Veni.202.107, and in part by Agence Nationale de la Recherche under the MELODIA project, grant number ANR-20-CE40-0013.

\section{Notation and definitions}\label{sec:notation}
All rings in this paper are commutative and unitary.
Let~$Z$ be a Dedekind domain with fraction field~$Q$.
For us, fields are not Dedekind domains.
Let~$V$ be a finite dimensional~$Q$-vector space.
A \emph{lattice} in~$V$ is a finitely-generated sub-$Z$-module of~$V$ which contains a~$Q$-basis of~$V$.
In particular, if~$L$ is such a lattice then~$LQ=V$.

Let~$K_1,\ldots,K_n$ be finite field extensions of~$Q$ and let~$K$ be the direct product
\[ K = K_1 \times \cdots \times K_n. \]
An \emph{order}~$R$ in~$K$ is a subring of~$K$ which is also a lattice in~$K$.
A \emph{fractional~$R$-ideal} is a finitely generated sub-$R$-module of~$K$ which is also a lattice in~$K$, or, equivalently, which contains an invertible element of~$K$.
For example, if~$R$ is a domain then an ideal~$I$ of~$R$ is a fractional~$R$-ideal if and only if~$I$ is non-zero.
Two fractional~$R$-ideals~$I$ and~$J$ are called \emph{isomorphic} if they are so as~$R$-modules. 
This is equivalent to having an invertible element~$\lambda$ in~$K$ such that~$I=\lambda J$.
Observe that given two fractional~$R$-ideals~$I$ and~$J$, the sum~$I+J$, the product~$IJ$ and the colon
\[ (I:J) = \set{ x\in K\ :\ xJ\subseteq I} \]
are also fractional~$R$-ideals.
The \emph{multiplicator ring} of a fractional~$R$-ideal~$I$ is the order~$(I:I)$.
We say that a fractional~$R$-ideal~$I$ is \emph{invertible} if there exists a fractional~$R$-ideal~$J$ such that~$IJ=R$.
Note that if this is the case, then~$J$ equals~$(R:I)$.
We define the \emph{Picard group} of~$R$ as the abelian group~$\Pic(R)$ of isomorphism classes of invertible fractional~$R$-ideals with the operation of multiplication.

Since~$K$ is commutative, the integral closure~$\cO$ of~$Z$ (embedded diagonally) in~$K$ is a subring of~$K$.
More precisely, we have
\[ \cO = \cO_1\times\cdots\times \cO_n, \]
where~$\cO_i$ is the integral closure of~$Z$ in~$K_i$.
By the Krull-Akizuki Theorem, each~$\cO_i$ is a Dedekind domain.
Every element of any order is integral over~$Z$,
which implies that~$\cO$ contains every order in~$K$.

Define the \emph{conductor}~$\frf$ of an order~$R$ in~$\cO$ as
\[ \frf = (R:\cO) = \set{ x\in K : x\cO \subseteq R }. \]
Note that since~$\frf$ is an~$\cO$-ideal, we have a decomposition
\[ \frf=\frf_1\oplus\cdots \oplus\frf_n, \]
where each~$\frf_i$ is an~$\cO_i$-ideal.
One can prove that~$\frf$ contains an invertible element of~$K$ if and only if~$\cO$ is finitely generated as a~$Z$-module.
If this is the case then~$\cO$ is an order, the \emph{maximal order} of~$K$.
We will assume that this is the case for the rest of the paper.
Moreover, for every fractional~$\cO$-ideal~$I$ we have a decomposition~$I=\oplus_{i=1}^n I_i$
for fractional~$\cO_i$-ideals~$I_i$, and hence also
\[ \Pic(\cO) = \Pic(\cO_1)\times \cdots \times\Pic(\cO_n). \]

\section{Isomorphism classes of lattices}
\label{sec:modules}
Let~$Z$,~$Q$,~$K=K_1\times \cdots \times K_n$ and~$\cO=\cO_1\times\cdots\times\cO_n$ be as in Section~\ref{sec:notation}.
We assume that~$\cO$ is finitely generated as a~$Z$-module.
Fix positive integers~$s_1,\ldots,s_n$ and consider the~$K$-module
\[ V = K_1^{s_1}\oplus \cdots \oplus K_n^{s_n}, \]
where the action of~$K$ is component-wise diagonal.

Let~$R$ be an order in~$K$.
Denote by~$\cL(R,V)$ the category of sub-$R$-modules of~$V$ which are also lattices, with~$R$-linear morphisms.
For every~$M$ in~$\cL(R,V)$, the extension~$M\cO$ is also a lattice in~$V$ and hence it belongs to~$\cL(\cO,V)$.
Pick a morphism~$\vphi:M\to N$ in~$\cL(R,V)$.
Since~$RQ=\cO Q=K$ and~$MQ=NQ=V$, the morphism~$\vphi$ extends uniquely to a~$K$-linear endomorphism of~$V$, which in turns restricts to a unique morphism~$M\cO\to N\cO$ in~$\cL(\cO,V)$.
We denote the induced morphisms also by~$\vphi$. 

The following is a restatement of Steinitz Theory; see~\cite{Steinitz12}.
\begin{prop}\label{prop:Steinitztheory}
    Let~$M$ be in~$\cL(\cO,V)$.
    Then there are fractional~$\cO_i$-ideals~$I_i$ and there exists an~$\cO$-linear isomorphism
    \[ M\simeq
    \bigoplus_{i=1}^n \left(\cO_i^{\oplus(s_i-1)}\oplus I_i\right).
    \]
    Moreover, the isomorphism class of~$M$ is uniquely determined by the integers~$s_i$ and the isomorphism class of the fractional~$\cO$-ideal~$I=I_1\oplus \cdots \oplus I_n$.
\end{prop}
\begin{proof}
    Since~$\cO=\cO_1\times\cdots\times\cO_n$ is a direct product of rings, we have a decomposition
    \[ M = M_1\oplus\cdots\oplus M_n, \]
    where each~$M_i$ is in~$\cL(\cO_i,K_i^{s_i})$.
    By the structure theorem for finitely generated modules over a Dedekind domain, see~\cite{Steinitz12}, there exists a fractional~$\cO_i$-ideal~$I_i$ such that
    \[ M_i \simeq \cO_i^{\oplus(s_i-1)}\oplus I_i. \]
    Furthermore, the isomorphism class of each~$M_i$ is uniquely determined by the rank~$s_i$ and the isomorphism class of~$I_i$.
    The same holds true after taking the direct sum.
\end{proof}

We now consider~$\cL(R,V)$ for a general order~$R$.
Let~$\frf=\oplus_{i=1}^n\frf_i$ be the conductor of~$R$ in~$\cO$.
\begin{thm}\label{thm:bwcondandmaxord}
    Let~$M$ be in~$\cL(R,V)$.
    Then there exist an~$M'$ in~$\cL(R,V)$, and fractional~$\cO_i$-ideals~$I_i$ such that
    \begin{enumerate}[(i)]
        \item~$M'\simeq M$ as an~$R$-module.
        \item~$M'\cO = \bigoplus_{i=1}^n \left(\cO_i^{\oplus(s_i-1)}\oplus I_i\right)$.
        \item~$\bigoplus_{i=1}^n \left(\frf_i^{\oplus(s_i-1)}\oplus \frf_iI_i\right) \subseteq M' \subseteq
        \bigoplus_{i=1}^n \left(\cO_i^{\oplus(s_i-1)}\oplus I_i\right)$.
    \end{enumerate}
\end{thm}
\begin{proof}
    By Proposition~\ref{prop:Steinitztheory}, there exists an~$\cO$-linear isomorphism
    \[ \vphi:M\cO \overset{\sim}{\longrightarrow} 
    \bigoplus_{i=1}^n \left(\cO_i^{\oplus(s_i-1)}\oplus I_i\right), \]
    for some fractional~$\cO_i$-ideals~$I_i$.
    Define~$M'=\vphi(M)$.
    Observe that~$M'$ is in~$\cL(R,V)$ and~$M\simeq M'$ by construction.
    We have
    \[  M' \subseteq M'\cO = \vphi(M)\cO = \vphi(M\cO) = \bigoplus_{i=1}^n \left(\cO_i^{\oplus(s_i-1)}\oplus I_i\right).\]
    For the other inclusion, note that
    \[
        \bigoplus_{i=1}^n \left(\frf_i^{\oplus(s_i-1)}\oplus \frf_iI_i\right) =
        \frf \left(\bigoplus_{i=1}^n \left(\cO_i^{\oplus(s_i-1)}\oplus I_i\right)\right) =
        \frf M'\cO = 
        \frf M' \subseteq M',
    \]
    where the last equality holds because~$\frf\cO=\frf$, and the inclusion follows from~$\frf \subseteq R$ and~$M'R=M'$.
\end{proof}

In the rest of the section we describe how to turn Theorem~\ref{thm:bwcondandmaxord} into an algorithm to compute representatives of the isomorphism classes of modules in~$\cL(R,V)$. 
We will need some additional assumptions on~$R$ and~$\cO$, which will be discussed in Remark~\ref{rmk:hyp}. We need to assume that:
\begin{enumerate}[(A)]
    \item \label{hyp:Zlattices} We have algorithms for working with $Z$-lattices in $V$ and $K$ and for working with fractional $R$-ideals for an arbitrary order in $K$.
    \item \label{hyp:maximalorder} We have an algorithm to compute the maximal order $\cO$ of $K$ and the conductor $(R:\cO)$ of an arbitrary order $R$ in $K$.
    \item \label{hyp:picardgroups} For each~$i$, the group~$\Pic(\cO_i)$ is finite, and we have an algorithm \texttt{PicardGroup} to compute it.
    \item \label{hyp:finmanysubmod} For each fractional~$\cO$-ideal~$I=I_1\oplus\cdots\oplus I_n$, the quotient
    \[ \cQ(I) = \dfrac{\cO_1^{\oplus(s_1-1)}\oplus I_1\oplus\cdots\oplus\cO_n^{\oplus(s_n-1)}\oplus I_n}{\frf_1^{\oplus(s_1-1)}\oplus \frf_1I_1 \oplus \cdots \oplus\frf_n^{\oplus(s_n-1)}\oplus \frf_nI_n} \]
    has finitely many sub-$R$-modules, and we have an algorithm called \texttt{SubModules} to list them all.
    \item \label{hyp:isiso} We have an algorithm \texttt{IsIsomorphic} that returns whether~$M$ and~$N$ in~$\cL(R,V)$ are isomorphic.
\end{enumerate}

\begin{remark}\label{rmk:hyp}
    We will focus on the case when $Q$ is a global field, that is, a number field or a finite extension of a function field~$k(T)$ where~$k$ is a finite field and~$T$ is an indeterminate.
    Under this assumption, the number of isomorphism classes of $\cL(R,V)$ is finite by the Jordan-Zassenhaus Theorem \cite[Thm.~26.4]{Reiner03}, which holds in much greater generality than the case we are considering.
    For example, it holds also for non-commutative orders.
    Nevertheless, Theorem~\ref{thm:bwcondandmaxord} could be interpreted as an effective version of the Jordan-Zassenhaus Theorem in the special case of a commutative $Z$-order $R$ in an \'etale $Q$-algebra $K$ with $Q$ a global field.
    
    If $Q$~is a global field, hypothesis \ref{hyp:Zlattices} is satisfied by using algorithms based on (pseudo) hermite normal form or linear algebra over finite fields. See for example \cite{cohen93} and \cite{cohenadv00}.
    In the same context also hypothesis \ref{hyp:maximalorder} is satisfied: see for example \cite{Bauch16} for the computation of the maximal order and see \cite[Sec.~6]{klupau05} for the conductor.

    Under the running assumption that~$Q$ is a global field, $\Pic(\cO)$ is a finite abelian group.
    If $Q$ is a number field, the problem of computing each~$\Pic(\cO_i)$ is classical, see~\cite{BuchWill89}.
    For finite extensions of function field, see~\cite{Hess99}.
    Hence, in both cases, hypothesis~\ref{hyp:picardgroups} is satisfied.
    
    Assume again that~$Q$ is a global field.
    Then the quotient~$\cQ(I)$ defined in~\ref{hyp:finmanysubmod} is a finite abelian group.
    We want to list all sub-$R$-modules~$N$ of~$\cQ(I)$ with trivial extension, that is, such that~$N\cO = \cQ(I)$.
    These modules are in bijection with the modules~$M'$ from Theorem~\ref{thm:bwcondandmaxord}.
    We will produce them by recursively computing the sub-$R$-modules of~$\cQ(I)$ which are maximal with respect to inclusion, as we now explain.
    This procedure is an adaptation of~\cite[Sec.~5.2]{FiekHofSirc19}.
    Let~$N$ be a such a maximal sub-$R$-module.
    Then there exists a rational prime~$p$ such that~$p\cQ(I)\subseteq N$.
    Hence~$N$ can be identified with a sub-$\F_p$-vector space of~$\cQ(I)/p\cQ(I)$.
    Now, one can use the \texttt{MEATAXE} algorithm, see~\cite{Par84} and~\cite[Sec.~7.4]{HoltEickOBrien05}, to enumerate all the maximal sub-$R$-modules of~$\cQ(I)/p\cQ(I)$, which are closed under the induced action of~$R$.
    From this list we need to keep only the ones with trivial extension: indeed if~$N\cO \neq \cQ(I)$ then all the sub-$R$-modules of~$N$ will not have trivial extension as well.
    Now we repeat the process with~$N$ instead of~$\cQ(I)$, and so on recursively until we have all sub-$R$-modules of~$\cQ(I)$ with trivial extension.

    Finally, in~\cite{2022BleyHofmannJohnston}, the authors describe an algorithm \texttt{IsIsomorphic} that in particular works for orders in \'etale algebras over~$\Q$.
    Hence, such orders satisfy also assumption~\ref{hyp:isiso}.
    To the best of our knowledge, there is no known analogous algorithm when~$Q$ is an extension of a function field.

    We conclude that orders in \'etale algebras over~$\Q$ satisfy all hypotheses \ref{hyp:Zlattices}, \ref{hyp:maximalorder}, \ref{hyp:picardgroups}, \ref{hyp:finmanysubmod} and~\ref{hyp:isiso}.
    Note also that if~$Q$ is a global field, then~$Z$ is Japanese, which means that~$\cO$ is finitely generated as a~$Z$-module, as required at the beginning of the section.
\end{remark}

We introduce now some notation that we will use throughout the rest of this section.
For any fractional~$\cO$-ideal~$I=\oplus_i I_i$ define the quotient
\[ \cQ(I) = \dfrac{\cO_1^{\oplus(s_1-1)}\oplus I_1\oplus\cdots\oplus\cO_n^{\oplus(s_n-1)}\oplus I_n}{\frf_1^{\oplus(s_1-1)}\oplus \frf_1I_1 \oplus \cdots \oplus\frf_n^{\oplus(s_n-1)}\oplus \frf_nI_n}, \]
as done before in~\ref{hyp:finmanysubmod}.
Denote by~$q_I$ the quotient map onto~$\cQ(I)$.
Define
\[ \widetilde\cM_I=\set{ \text{sub-$R$-module~$\widetilde N$ of~$\cQ(I)$}: \widetilde N\cO=\cQ(I) },\]
and
\[ \cM_I=\set{ q_I^{-1}(\widetilde N)\in \cL(R,V): \widetilde N \in \widetilde\cM_I }.\]

The following lemma will be used in Theorem~\ref{thm:algcorrect} to prove the correctness of Algorithm~\ref{alg:isom_classes}.
\begin{lemma}\label{lemma:algcorrect}\
    \begin{enumerate}[(i)]
        \item \label{lemma:algcorrect:isom}
        There is an~$R$-linear isomorphism~$\widetilde\psi: \cQ(\cO) \to \cQ(I)$.
        \item \label{lemma:algcorrect:tildepsi} 
        The isomorphism~$\widetilde\psi$ induces a bijection~$\psi: \cM_\cO \to \cM_I$ defined by 
       ~$\psi(M) = q_I^{-1}(\widetilde\psi(q_\cO(M)))$.
        \item \label{lemma:algcorrect:psi} 
        The bijection~$\psi$ induces a bijection between the sets of~$R$-linear isomorphism classes in~$\cM_\cO$ and in~$\cM_I$.     
    \end{enumerate}
\end{lemma}
\begin{proof}
    By replacing~$I$ with an isomorphic fractional~$\cO$-ideal, we can assume that~$I$ is coprime to the conductor~$\frf$, that is,~$\cO=\frf + I$; see~\cite[Cor.~1.2.11]{cohenadv00}.
    This implies that~$\frf_iI_i=\frf_i\cap I_i$ and~$\frf_i+I_i = \cO_i$, for every~$i$.
    Therefore, for every~$i$, we have the following isomorphism of~$R$-modules
    \[ \dfrac{I_i}{\frf_i I_i}=\dfrac{I_i}{\frf_i\cap I_i}\simeq \dfrac{I_i+\frf_i}{\frf_i} = \dfrac{\cO_i}{\frf_i}. \]
    By taking direct sums, we obtain the desired isomorphism~$\widetilde\psi:\cQ(\cO)\overset{\sim}{\to}\cQ(I)$, completing the proof of Part~\ref{lemma:algcorrect:isom}.
    Part~\ref{lemma:algcorrect:tildepsi} is an immediate consequence of Part~\ref{lemma:algcorrect:isom}.
    For Part~\ref{lemma:algcorrect:psi} we argue as follows.
    Pick~$M_1$ and~$M_2$ in~$\cM_\cO$ and let~$\vphi:M_1\overset\sim\to M_2$ be an~$R$-linear isomorphism.
    Since~$\cQ(\cO) = M_i\cO\otimes_\cO (\cO/\frf)$, for~$i=1,2$,
    we see that~$\vphi$ induces an automorphism of~$\cQ(\cO)$.
    Pushing this forward via~$\widetilde\psi$, we obtain an automorphism of~$\cQ(I)$ which then lifts to an isomorphism~$\psi(M_1)\simeq \psi(M_2)$.
    In an analogous manner, given isomorphic~$M_1$ and~$M_2$ in $\cM_I$, we obtain an isomorphism $\psi^{-1}(M_1)\simeq \psi^{-1}(M_2)$.
    Hence we obtain a bijection between the sets of~$R$-linear isomorphism classes in~$\cM_\cO$ and in~$\cM_I$, as required.
\end{proof}

\begin{alg}\label{alg:isom_classes}
    Assume that~\ref{hyp:Zlattices},~\ref{hyp:maximalorder},\ref{hyp:picardgroups},~\ref{hyp:finmanysubmod} and~\ref{hyp:isiso} hold.
    The following steps will return a list~$\cL_{\mathrm{out}}$ of representatives of the isomorphism classes of~$\cL(R,V)$.
    \begin{enumerate}[(1)]
        \item \label{step:conductor} Compute the maximal order $\cO$ of $K$ and the conductor~$\frf=\frf_1\oplus\cdots \oplus\frf_n$ of~$R$ in~$\cO$.
        \item \label{step:class_gp} Use \normalfont\texttt{PicardGroup} to compute representatives~$I^{(k)} = \oplus_{i=1}^n I^{(k)}_i~$ of
        \[ \Pic(\cO)=\bigoplus_{i=1}^n \Pic(\cO_i). \]
        \item \label{step:quotient} Form the quotient
        \[ \cQ(\cO) = \dfrac{\cO_1^{\oplus s_1}\oplus\cdots\oplus\cO_n^{\oplus s_n}}{\frf_1^{\oplus s_1}\oplus\cdots\oplus\frf_n^{\oplus s_n}}, \]
        and denote by~$q_\cO$ the natural quotient map.
        \item \label{step:enum} Use \normalfont\texttt{SubModules} to produce a list~$\widetilde\cM_\cO$ of all the finitely many sub-$R$-modules~$\widetilde M$ of~$\cQ(\cO)$ such that~$\widetilde M\cO=\cQ(\cO)$.
        \item Initialize an empty list~$\cL$.
        \item \label{step:sieving} For each~$\widetilde M$ in~$\widetilde\cM_\cO$ do:
            \begin{enumerate}[(a)]
                \item Compute~$M=q_\cO^{-1}(\widetilde M)$.
                \item If 
                \normalfont\texttt{IsIsomorphic} returns that there is no module~$M'$ in~$\cL$ which is isomorphic to~$M$ then append~$M$ to~$\cL$.
            \end{enumerate}
        \item Initialize an empty output list~$\cL_{\mathrm{out}}$.
        \item \label{step:foreachk} For each~$k$ do:
            \begin{enumerate}[(a)]
                \item Initialize an empty list~$\cL^{(k)}$.
                \item Compute~$\widetilde\psi:\cQ(\cO)\overset{\sim}{\to} \cQ(I^{(k)})$ as in Lemma~\ref{lemma:algcorrect}.\ref{lemma:algcorrect:isom}.
                \item For each~$M$ in~$\cL$ do:
                    \begin{enumerate}
                        \item Compute~$M'=q_I^{-1}(\widetilde\psi(q_\cO(M)))$.
                        \item Append~$M'$ to~$\cL^{(k)}$.
                    \end{enumerate}
                \item Concatenate~$\cL^{(k)}$ to~$\cL_{\mathrm{out}}$.
            \end{enumerate}
        \item Return~$\cL_{\mathrm{out}}$.
    \end{enumerate}
\end{alg}

\begin{thm}\label{thm:algcorrect}
    Algorithm~\ref{alg:isom_classes} returns a minimal set of representatives of the isomorphism classes in~$\cL(R,V)$.
\end{thm}
\begin{proof}
    By Theorem~\ref{thm:bwcondandmaxord} for every~$M$ in~$\cL(R,V)$ there exist an index~$k$ and an element~$M'$ of~$\cM_{I^{(k)}}$ such that~$M\simeq M'$.
    By Lemma~\ref{lemma:algcorrect}.\ref{lemma:algcorrect:psi}, there exists~$M''$ in~$\cL^{(k)}$ such that~$M\simeq M''$.
    In other words, the list~$\cL_{\mathrm{out}}$ contains representatives of all isomorphism classes of~$\cL(R,V)$.
    We need to show that there are no repetitions, that is, that the elements of~$\cL_{\mathrm{out}}$ are pairwise non-isomorphic.
    Again by Lemma~\ref{lemma:algcorrect}.\ref{lemma:algcorrect:psi}, for each index~$k$, the elements of~$\cL^{(k)}$ are pairwise non-isomorphic.
    We are left to show that given two distinct indices~$k$ and~$k'$, and modules~$M\in \cL^{(k)}$ and~$M'\in \cL^{(k')}$, there cannot be an isomorphism~$M\simeq M'$. 
    Indeed, if this were the case then we would have an isomorphism~$M\cO\simeq M'\cO$, contradicting Proposition~\ref{prop:Steinitztheory}.
\end{proof}
\begin{cor}\label{cor:divbyclassnum}
    The number of isomorphism classes in~$\cL(R,V)$ is divisible by the size of~$\Pic(\cO)$.
\end{cor}
\begin{proof}
    This is a consequence of the fact that the lists~$\cL^{(k)}$ defined in Step~\ref{step:foreachk} of Algorithm~\ref{alg:isom_classes} all have size equal to the size of the list~$\cL$ built in Step~\ref{step:sieving}.
\end{proof}

\begin{remark}\label{rmk:is_isom}
    In Step~\ref{step:sieving} of Algorithm~\ref{alg:isom_classes}, we use as a black-box the algorithm \texttt{IsIsomorphic} to test whether two $R$-modules $M$ and $M'$  are isomorphic, where $M$ is fixed while $M'$ loops over all elements of $\cL$.
    As pointed out in Remark~\ref{rmk:hyp}, in the case of orders in \'etale algebras over $\Q$, we can use the algorithm provided by \cite{2022BleyHofmannJohnston}, which we now briefly review.
    Set $A=\Hom_{K}(V,V)$, $X=\Hom_R(M,M')$ and $\Lambda=\End_R(M')$.
    By~\cite[Prop~3.1]{2022BleyHofmannJohnston}, we have that $M$ and $M'$ are isomorphic if and only if the $\Lambda$-lattice $X$ is free of rank $1$, and every (any) free generator of $X$ over $\Lambda$ is an isomorphism.
    The number of the Steps below refers to \cite[Alg.~8.3]{2022BleyHofmannJohnston}, which returns whether $X$ is free of rank $1$ over $\Lambda$ and, if so, a generator.
    All outputs of Steps $(1)$-$(4)$, $(8)$-$(9)$ should be cached, since they depend only on $M$.
    Moreover, following the proof of~\cite[Thm.~8.4]{2022BleyHofmannJohnston},
    Step $(5)$ is probabilistic polynomial-time reducible to one call of $\texttt{IsPrincipal}$ for each $1\leq i\leq n$.
    Hence, one should also cache $\Pic(\cO_{K_i})$, and also $\cO_{K_i}^\times$ since it is used in Step $(10)$.
    Step $(6)$ is probabilistic polynomial-time reducible to computing a factorization of an integer, which should also be stored.
\end{remark}

\begin{remark}\label{running_time}
    We discuss here the running time of Algorithm~\ref{alg:isom_classes} under the assumption that $R$ is a $\Z$-order in an \'etale algebra over $\Q$.
    In this case the required operations with $\Z$-lattices and fractional ideals (cf.~hypothesis \ref{hyp:Zlattices}) can be performed in polynomial time.
    Steps~$(7)$ and~$(9)$ are trivial.
    Steps $(1)$ and $(3)$ have polynomial running time in the size of the input once the computation of the maximal order $\cO$ is completed which requires to know the prime factors of the discriminant of a polynomial defining $K$ over $\Q$.
    Step $(2)$ can be achieved in heuristic sub-exponential running time using for example \cite{CohDiazOlivier97}.
    Step~$(8)$ requires the computations of a representative coprime to the conductor which is probabilistic polynomial time; see \cite[Cor.~A.2]{2022BleyHofmannJohnston}.
    In Step $(6)$, we have to run \texttt{IsIsomorphic} a number of times which is bounded by above by the size of the output of Step $(4)$ and the size of the list $\cL$.
    The algorithm provided by \cite{2022BleyHofmannJohnston} reduces in probabilistic polynomial time to well known algorithm in number theory (like \texttt{IsPrincipal} and \texttt{UnitGroup}).
    It follows that the running time of Algorithm~\ref{alg:isom_classes} depends on the size of output of Step $(4)$, the list $\widetilde\cM_\cO$, which is produced by \texttt{SubModules} using, for example, the ideas described in Remark~\ref{rmk:hyp}.
    In general, it is not easy to give an upper bound for the size of this list, but a lower bound is computed as follows.
    Let $\p$ be a maximal ideal of $R$ above the conductor $\frf$ of $R$.
    Write
    $\p \cO = \mathfrak{P}_1 \oplus \cdots \oplus \mathfrak{P}_n $, where $\mathfrak{P}_i$ is an ideal of $\cO_i$, for each $i$.
    Then the number of sub-$R$-modules of $\cQ(\cO)$ is bounded from below from the number $C_\p$ of sub-$R/\p$-vector subspaces of the $R/\p$-vector space 
    \[ \frac{\cO_1^{\oplus s_1}\oplus \cdots \oplus \cO_n^{\oplus s_n}}{ \mathfrak{P}_1^{\oplus s_1} \oplus \cdots \oplus \mathfrak{P}_n^{\oplus s_n} }. \]
    Then taking the maximum of $C_\p$ over all maximal $R$-ideals above $\frf$ gives a lower bound.
\end{remark}

\begin{remark}\label{rmk:labelling}
    When~$R$ is an order in an \'etale algebra over~$\Q$, the process described in Algorithm~\ref{alg:isom_classes} can be used to label the isomorphism classes of lattices in~$\cL(R,V)$, as we now describe.

    The representatives of each~$\Pic(\cO_i)$ can be ordered in a deterministic way by \cite{2020CrePageSuth_arXiv}.
    Taking direct sums, we obtain an induced ordering on~$\Pic(\cO)$, and hence we can order the various quotients~$\cQ(I^{(k)})$.

    To conclude it suffices to sort the sub-$R$-modules~$\widetilde M$ of~$\cQ(\cO)$, that is, the elements of the list~$\widetilde\cM_\cO$.
    This can be done, for example, as follows.
    Let~$q_\cO$ be the quotient map onto~$\cQ(\cO)$.
    For each~$\widetilde M$, compute the Hermite Normal Form of the matrix representing a~$\Z$-basis of~$M_0:=q_\cO^{-1}(\widetilde M)$ with respect to a fixed~$\Q$-basis of~$V$.
    Then one simply sorts these matrices to obtain the desired result.
    In fact, one can use this method together with invariants that take into account the~$R$-module structure of~$M_0$, like, for example, the~$R/\p$-dimensions of~$M_0\otimes_R R/\p$ where~$\p$ runs over the finitely many primes~$\p$ of~$R$ above the conductor~$\frf=(R:\cO)$.
\end{remark}

\section{Conjugacy classes of semisimple integral matrices}
\label{sec:matrices}
Let~$A$ and~$B$ be two square matrices with integer coefficients, both of dimension~$r$.
Recall that~$A$ and~$B$ are~$\Z$-conjugate if there exists a matrix~$P$ in~$\GL_r(\Z)$ such that~$PA=BP$, in which case we will write~$A \sim_\Z B$.
If~$A$ and~$B$ are~$\Z$-conjugate then they have the same minimal and characteristic polynomials.

Let~$m$ be a squarefree polynomial in~$\Z[x]$, that is, such that the factors appearing in its irreducible factorization
\[ m = m_1\cdots m_n \]
generate pairwise coprime ideals in $\Q[x]$.
Fix positive integers~$s_1,\ldots,s_n$.
Define
\[ h = m_1^{s_1}\cdots m_n^{s_n}. \]
Denote by~$\Mat_{m,h}$ the set of integral square matrices with minimal polynomial~$m$ and characteristic polynomial~$h$.
Since~$m$ is squarefree, these matrices are semisimple.
Consider the \'etale algebra~$K= \Q[x]/m$, and the order~$R=\Z[\pi]$, where~$\pi$ denotes the class of the variable~$x$ in~$K$.
Consider
\[ V = K_1^{s_1}\oplus \cdots \oplus K_n^{s_n}. \]
As before,~$\cL(R,V)$ denotes the category of~$\Z$-lattices in~$V$ which are~$R$-modules.
Pick~$M$ in~$\cL(R,V)$ and choose a~$\Z$-basis~$\cB$ of~$M$. 
Define~$A_{M,\cB}$ as the matrix that represents multiplication by~$\pi$ on~$M$ with respect to the basis~$\cB$.
Since~$M$ is an~$R$-module, the matrix~$A_{M,\cB}$ has integer entries.
Denote by~$\Psi$ the function that associates the pair~$(M,\cB)$ to the matrix~$A_{M,\cB}$.
In previous work, we proved the following theorem, which is a generalization of the Latimer-MacDuffee Theorem~\cite{LaClMD33}.
\begin{thm}{\cite[Thm.~8.1]{MarICM18}}\label{thm:LM}
    The function~$\Psi$ induces a bijection between the isomorphism classes in~$\cL(R,V)$ and~$\Mat_{m,h}/\sim_\Z$.
\end{thm}

Combining Theorem~\ref{thm:LM} with Algorithm~\ref{alg:isom_classes} we obtain the following theorem.
\begin{thm}\label{thm:conjclasses} 
    Algorithm~\ref{alg:isom_classes} allows us to compute a minimal set of representatives of~$\Mat_{m,h}/\sim_\Z$.
\end{thm}

\begin{remark}\label{rmk:matrix_decomposition}
    As mentioned above, the matrices in~$\Mat_{m,h}$ are semisimple.
    It is known that, in general, an integral square matrix~$A$ can be written uniquely as~$A=S+N$, where~$S$ is semisimple,~$N$ is nilpotent and~$SN=NS$. 
    If~$A'=S'+N'$ is another matrix, with analogous decomposition, then there exists an invertible integral matrix~$P$ such that~$PA=A'P$ if and only if~$PS=S'P$ and there exists an invertible integral matrix~$T$ in the stabilizer of~$S'$ such that~$TPN=N'TP$.
    This seems to suggest that in order to generalize Theorem~\ref{thm:LM}, one first needs to find a method to compute conjugacy classes of nilpotent matrices with prescribed minimal and characteristic polynomials, where the conjugation is realized only by matrices in a subgroup of~$\GL_n(\Z)$.
    To the best of our knowledge, this problem has not been solved yet.
\end{remark}

\begin{remark}
    Theorem~\ref{thm:LM} together with Remark~\ref{rmk:labelling} gives a method to label the representatives of the conjugacy classes in~$\Mat_{m,h}$.
\end{remark}

\section{Isomorphism classes of abelian varieties over finite fields}
\label{sec:AV}
In this section we describe how to use Algorithm~\ref{alg:isom_classes} to compute the isomorphism classes of abelian varieties over a finite field~$\F_q$ belonging to isogeny classes satisfying certain hypotheses.
Here, by isogeny classes, we mean~$\F_q$-isogeny classes.
Recall that by Honda-Tate theory, see~\cite{Tate66}, \cite{Honda68} and~\cite{Tate71}, such an isogeny class is uniquely determined by the characteristic polynomial~$h$ of Frobenius of any abelian variety in the isogeny class.
We will denote the isogeny class by~$\cI_h$, and turn it into a category by considering~$\F_q$-morphisms between the objects.
The polynomial~$h$ is in~$\Z[x]$, has degree~$2g$, where~$g$ is the dimension of any abelian variety in~$\cI_h$, and all its complex roots have norm~$\sqrt{q}$.
Recall also that~$\cI_h$ is ordinary if the coefficient of~$x^g$ in~$h$ is coprime to~$q$.

Consider the factorization
\[ h=m_1^{s_1} \cdots m_n^{s_n} \]
into irreducible factors, with the $m_i$ generating pairwise coprime ideals in~$\Q[x]$.
Put~$m=m_1\cdots m_n$.
Define~$K_i=\Q[x]/m_i$ for each~$i$, and denote by~$\pi$ the class of~$x$ in~$K = K_1 \times \cdots \times K_n$.
Consider the order~$R=\Z[\pi,q/\pi]$ in~$K$.
Finally, set
\[ V = K_1^{s_1}\oplus \cdots \oplus K_n^{s_n}. \]
As before,~$\cL(R,V)$ denotes the category of~$\Z$-lattices in~$V$ which are~$R$-modules, with~$R$-linear morphisms.

\begin{thm}\label{thm:eqcatAV}\
    \begin{enumerate}[(i)]
        \item Assume that~$h$ is ordinary.
        Then there is a categorical equivalence between~$\cI_h$ and~$\cL(R,V)$.
        \item Assume that~$q$ is prime and that~$h$ does not have real roots.
        Then there is a categorical anti-equivalence between~$\cI_h$ and~$\cL(R,V)$.
    \end{enumerate}
    In both cases, we obtain a bijection between the isomorphism classes of abelian varieties in~$\cI_h$ and the isomorphism classes in~$\cL(R,V)$.
\end{thm}
\begin{proof}
    See~\cite[Theorem 4.1.(a)]{MarBassPow} for a complete proof in the case where~$h$ is a power of~$m$.
    We review the main steps here for completeness.
    We first handle the ordinary case.
    By Deligne~\cite[Sec.~7, Th\'eor\`eme]{Del69}, there is an equivalence between~$\cI_h$ and the category of pairs~$(T,F)$ where~$T$ is a free finitely generated~$\Z$-module and~$F$ is an endomorphism of~$T$ with characteristic polynomial~$h$, which is semisimple, and there is a~$\Z$-linear endomorphism~$V$ of~$T$ such that~$V\circ F$ and~$F\circ V$ are both the multiplication-by-$q$ map on~$T$.
    A morphism between two pairs~$(T,F)$ and~$(T',F')$ is a~$\Z$-linear morphism~$\vphi:T\to T'$ such that~$F'\circ \vphi = \vphi \circ F$.
    
    Fix a pair~$(T,F)$ as above.
    Since~$F$ is semisimple, its minimal polynomial is~$m$.
    Hence we can identify~$F$ with~$\pi$ and~$V$ with~$q/\pi$.
    This identification induces an~$R$-module structure on~$T$, and induces a~$K$-linear isomorphism~$T\otimes_\Z \Q \simeq V$.
    Via this isomorphism, we can identify~$T$ with an element of~$\cL(R,V)$.
    Therefore we obtain an equivalence between~$\cI_h$ and~$\cL(R,V)$.

    Applying the same argument to the result of Centeleghe-Stix~\cite[Thm.~1]{CentelegheStix15}, with arrows reversed, induces the anti-equivalence in the case where~$q$ is prime and~$h$ has no real roots.
\end{proof}

Combining Theorem~\ref{thm:eqcatAV} with Algorithm~\ref{alg:isom_classes} we obtain the following result.
\begin{thm}\label{thm:isomclassesAV}
    Assume that~$h$ is ordinary, or that~$q$ is prime and~$h$ has no real roots. 
    Algorithm~\ref{alg:isom_classes} allows us to compute all the representatives of the isomorphism classes of abelian varieties in~$\cI_h$.
\end{thm}

\begin{remark}
    Theorem~\ref{thm:isomclassesAV} together with Remark~\ref{rmk:labelling} gives a method to label abelian varieties in the isogeny class~$\cI_h$. 
\end{remark}

\section{Examples}\label{sec:example}
Algorithm~\ref{alg:isom_classes} has been implemented in \texttt{MAGMA}~\cite{Magma}.
Such implementation is available at 
\url{https://github.com/stmar89/AlgEt}\footnote{at the moment of submission, the most recent commit is 4c22349}, together with the code to reproduce the examples (see the webpage of the author for a precise link).
Even if it is possible run the code using only \texttt{MAGMA}, it is much faster to  combine it with the implementation of \texttt{IsIsomorphic} from~\cite{2022BleyHofmannJohnston} which is included in the \texttt{julia} package \texttt{Nemo/Hecke}~\cite{nemo}.
Detailed instructions are included in the code.
\begin{example}\label{ex:AV}
    Consider the polynomials 
    \[ m_1=x^2 - x + 3 \quad \text{and}\quad  m_2=x^2 + x + 3. \]
    Put~$K_i=\Q[x]/m_i$ for~$i=1,2$ and let~$\pi_i$ be the class of~$x$ in~$K_i$.
    Consider the isogeny classes~$\cI_{m_1}$ and~$\cI_{m_2}$ of ordinary elliptic curves over~$\F_3$ determined by~$m_1$ and~$m_2$.
    Since both~$\cO_1=\Z[\pi_1]$ and~$\cO_2=\Z[\pi_2]$ are maximal and have trivial Picard group, we see that both isogeny classes contain a unique isomorphism class of elliptic curves.
    We will denote the representatives by~$E_1$ and~$E_2$, respectively.

    Now we consider the product
    \[ m=m_1m_2 = x^4 +5x^2 + 9.\]
    Put~$K=K_1\times K_2$ and denote by~$\pi$ the class of~$x$ in~$K$.
    Put~$R=\Z[\pi,3/\pi]$. 
    Using Algorithm~\ref{alg:isom_classes}, we see that there are~$4$ isomorphism classes in~$\cL(R,K)$.
    Note that such modules are fractional~$R$-ideals in~$K$, so in fact it is more efficient to use the specialized code from~\cite{MarICM18}.
    The only overorder of~$R$ is the maximal order~$\cO=\cO_1\times\cO_2$ of~$K$. 
    Among the~$4$ classes,~$3$ come from~$\Pic(R)$, and the last one from~$\Pic(\cO)$, which is trivial.
    These isomorphism classes correspond to~$4$ isomorphism classes of abelian varieties in~$\cI_{m}$ by Theorem~\ref{thm:eqcatAV}.
    More precisely, the isomorphism classes in~$\cI_{m}$ are represented by~$3$ abelian surfaces,~$A_1,A_2,A_3$, all with endomorphism ring isomorphic to~$R$, and then the class of~$E_1\times E_2$.
    Since the order~$R$ is not a direct product of orders from~$K_1$ and~$K_2$, we deduce that no surface~$A_i$ is isomorphic to a product of elliptic curves.

    Consider the isogeny class~$\cI_h$ of abelian threefolds determined by 
    \[ h=m_1^2m_2 = x^6 - x^5 + 8x^4 -5x^3 + 24x^2 - 9x - 27. \]
    Algorithm~\ref{alg:isom_classes} shows that there are~$2$ isomorphism classes in~$\cL(R,K_1^2\oplus K_2)$, represented by~$\cO_1\times R$ and~$\cO_1 \times \cO_1 \times \cO_2$.
    By Theorem~\ref{thm:eqcatAV}, these correspond to the isomorphism classes of abelian varieties in~$\cI_{h}$, represented by
    \[ A_1\times E_1\text{ and } E_1^2\times E_2. \]
    In particular, we have~$\F_3$-isomorphisms
    \[ A_1\times E_1 \simeq A_2\times E_1 \simeq A_3\times E_1. \]
\end{example}

\begin{example}\label{ex:matrices}
    In this example we consider again the polynomials
    \begin{gather*}
        m_1 =x^2 - x + 3, \\
        m_2 = x^2 + x + 3,\\
        m = m_1m_2 = x^4 +5x^2 + 9, \\
        h =m_1^2m_2 = x^6 - x^5 + 8x^4 -5x^3 + 24x^2 - 9x - 27.,
    \end{gather*}
    from Example~\ref{ex:AV}.
    For~$i=1,2$, define~$K_i=\Q[x]/m_i$.
    Put~$K=K_1\times K_2$ as before, but now we consider the order~$E=\Z[\pi]$, where~$\pi$ is the class of~$x$ in~$K$ (instead of~$R=\Z[\pi,3/\pi]$ as in Example~\ref{ex:AV}).

    We want to compute the representatives of~$\Z$-conjugacy classes of integral matrices having minimal polynomial~$m$ and characteristic polynomial~$h$, that is, in the notation from Section~\ref{sec:matrices}, we want to compute representatives of~$\Mat_{m,h}/\sim_{\Z}$.
    By Theorem~\ref{thm:LM} we have to compute the isomorphism classes of modules in~$\cL(E,K_1^2\oplus K_2)$.
    Algorithm~\ref{alg:isom_classes} returns that there are~$4$ classes, represented by the following matrices.
    \begin{align*}
        \begin{pmatrix}
            0 & 1 & 0 & 0 & 0 & 0 \\
            -3 & 1 & 0 & 0 & 0 & 0 \\
            0 & 0 & 0 & 1 & 0 & 0 \\
            0 & 0 & -3 & 1 & 0 & 0 \\
            0 & 0 & 0 & 0 & 0 & 1 \\
            0 & 0 & 0 & 0 & -3 & -1
            \end{pmatrix}, & &
            \begin{pmatrix}
            0 & 1 & 0 & 0 & 0 & 0 \\
            -3 & 1 & 0 & 0 & 0 & -1 \\
            0 & 0 & 0 & 1 & -2 & 0 \\
            0 & 0 & -3 & 1 & -1 & 2 \\
            0 & 0 & 0 & 0 & 0 & 1 \\
            0 & 0 & 0 & 0 & -3 & -1
            \end{pmatrix}\\
            \begin{pmatrix}
            0 & 1 & 0 & 0 & 0 & 2 \\
            -3 & 1 & 0 & 0 & 2 & 0 \\
            0 & 0 & 0 & 1 & 0 & 1 \\
            0 & 0 & -3 & 1 & 1 & 0 \\
            0 & 0 & 0 & 0 & 0 & 3 \\
            0 & 0 & 0 & 0 & -1 & -1
            \end{pmatrix}, & &
            \begin{pmatrix}
            0 & 1 & 0 & 0 & 0 & 1 \\
            -3 & 1 & 0 & 0 & 1 & -1 \\
            0 & 0 & 0 & 1 & -1 & 2 \\
            0 & 0 & -3 & 1 & 2 & 3 \\
            0 & 0 & 0 & 0 & 0 & 3 \\
            0 & 0 & 0 & 0 & -1 & -1
            \end{pmatrix}.
    \end{align*}
\end{example}

\bigskip\noindent\textbf{Competing interest statement}\\
The author has no competing interests to declare that are relevant to the content of this article.

\bigskip\noindent\textbf{Data availability statement}\\
No datasets were generated or analysed during the current study.

\bigskip\noindent\textbf{Published version}\\
This version of the article has been accepted for publication, after peer review but is not the Version of Record and does not reflect post-acceptance improvements, or any corrections. 
The Version of Record is available online at: \url{https://doi.org/10.1007/s40993-024-00584-9}

\bibliographystyle{alpha}
\renewcommand{\bibname}{References} 
\bibliography{references} 
\end{document}